\documentclass[pra,preprint]{revtex4}
\usepackage{graphicx}
\usepackage{bm}

\begin{document}
\author{V\'{\i}ctor Romero-Roch\'{\i}n\footnote{Telephone: + 52 (55) 5622 5096. {\it E-mail address:} romero@fisica.unam.mx}} 
\affiliation{Instituto de F\'{\i}sica, Universidad
Nacional Aut\'onoma de M\'exico, Apdo. Postal 20-364, M\'exico D.
F. 01000, Mexico. }

\title{A derivation of Benford's Law ... and a vindication of Newcomb}

\date{\today}
\begin{abstract}
We show how Benford's Law (BL) for first, second, ... , digits, emerges from the distribution of digits of numbers of the type $a^{R}$, with $a$ any real positive number and $R$ a set of real numbers uniformly distributed in an interval $\left[ \left. P\log_a 10, (P +1) \log_a 10 \right) \right.$ for any integer $P$. The result is shown to be number base and scale invariant. A rule based on the mantissas of the logarithms allows for a determination of whether a set of numbers obeys BL or not. We show that BL applies to numbers obtained from the {\it multiplication} or {\it division} of numbers drawn from any distribution. We also argue that (most of) the real-life sets that obey BL are because they are obtained from such basic arithmetic operations. We exhibit that all these arguments were discussed in the original paper by Simon Newcomb in 1881, where he presented Benford's Law.
\end{abstract}

\maketitle

\section{Introduction.}
 Benford's Law (BL) asserts that in certain sets of numbers, most of them of real-life origin, the first digit is distributed non-uniformly in the form
 \begin{equation}
 P_B^{(1)}(d) = \log_{10} \left( 1 + \frac{1}{d} \right) , \label{BL}
 \end{equation}
 where $d$ is the first digit of the number and $\log_{10}$ is the logarithm base 10. In other words, $P_B^{(1)}(d)$ is the fraction of the numbers with first digit $d$ in the given set. There are also forms of Benford's Law for second, third, etc., digits, namely $P_B^{(n)}(d)$. Table 1 shows the values of $P_B^{(1)}(d)$ for $d = 1, 2, \dots , 9$. 
  
 \begin{table}[htdp]
\begin{tabular}{|c|c|c|c|c|c|c|c|c|c|}
\hline \hline
$d$ & 1 & 2 & 3 & 4 & 5 & 6 & 7 & 8 & 9 \\ \hline

$P_B^{(1)}$ & $\> 0.3010 \>$& $\>  0.1761 \>$& $\>  0.1249 \>$& $\>  0.0969 \>$ & $\> 0.0792 \>$& $\> 0.0669 \>$ & $\> 0.0580 \>$ &$\>  0.0512 \>$& $\> 0.0458\>$ \\
\hline \hline
\end{tabular}
\caption{First digit Benford law.}
\end{table}

 BL has been found to be obeyed quite well in a variety of situations, many of them checked by Franck Benford himself\cite{Benford}. These sets include population census, stock markets indeces, utilities bills, tax returns, areas of rivers, physical and mathematical constants, and molecular weights, among others\cite{Benford,Fewster}. At first sight, the Law is certainly baffling and counterintuitive\cite{Lines} since one's naive intuition is that digits of numbers should be uniformly or randomly distributed. Although Franck Benford has been credited with the law for his work of 1938\cite{Benford}, the law was originally discovered by the astronomer Simon Newcomb in 1881\cite{Newcomb} as a follow up of the observation that the pages of tables of logarithms in his university library were worn out following BL, as given by equation (\ref{BL}). What is rarely told is that Newcomb {\it derived} Benford's Law. His demonstration for us may now look obscure, and probably just sketchy, because he used arguments that were not so difficult to those familiar with concepts of log tables ... and certainly we are not. We shall advance a plausible explanation of Newcomb's observation of the worn pages of the log tables and argue why many sets of real-life origin also obey BL; alas, this argument was also used by Newcomb.

We shall first prove a general result that appears to be known already\cite{Fewster,Hill,Pietronero,Raimi,others}, although to the best of our knowledge it has not been shown explicitely in the form here presented; we shall see that yields exactly all digit's Benford's distributions, allowing also for concluding that it is scale and number base invariant\cite{Hill}. We demonstrate that if $R$ is a set of real numbers uniformly distributed, then, the distributions of digits of $a^R$ obey BL for any real positive number $a$. Then, we discuss the main result of Newcomb, namely, the fact that a given set of numbers obeys BL if  the mantissas of their logarithms are uniformly distributed. We then analyze two main type of sequences of numbers that obey BL, those that are obtained from multiplication of numbers drawn from any distribution and those that are part of a geometric progression of numbers uniformly distributed in an arbitrary interval.

\section{A general result concerning Benford's Law.}

Let $\{R_1, R_2, \dots , R_N\}$ be a sequence of real numbers drawn from a uniform distribution in the interval $R_i \in \left[ \left. P \log_a 10,  (P +1) \log_a 10 \right) \right.$, with $P$ any integer. Then, the first, second, ..., digit distributions of the sequence $\{a^{R_1}, a^{R_2}, \dots , a^{R_N}\}$, with $a$ any real positive number, approaches Benford's Law, Eq.(\ref{BL}) and its generalizations, as $N \to \infty$ .

\begin{figure}[ htbp!] 
  \centering
\includegraphics[width=0.5\textwidth,keepaspectratio]{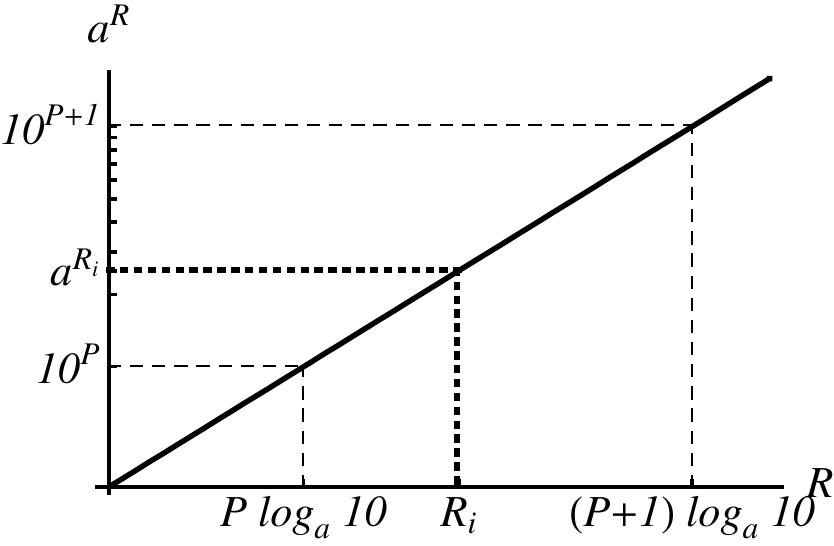}
\caption{$a^R$ vs $R$ in semi-log scale. The dotted line shows an example of an arbitrary point $R_i$ in the chosen interval, such that the first digit of $a^{R_i}$ is 2 because it falls in bin $b_2$.}
\label{GR-2}
\end{figure}

 Let us look first at the first digit distribution. In Fig. \ref{GR-2} we plot $a^R$ vs $R$ in a semi-log (base $a$) scale. In this graph, $a^R$ vs $R$ appears as a straight line. Now take in the $R$-axis the sequence $\{R_1, R_2, \dots , R_N\}$ within the interval $\left[ \left. P \log_a 10, (P +1) \log_a 10 \right) \right.$. Take any number of the sequence, say $R_i$. Then, in the logarithmic scale, $\log_a a^{R_i}$ must lie within any of the following ``bins": $b_1$, the interval $[\log_{a}\left(1\times 10^P\right), \log_{a} \left(2\times 10^P\right) )$; or $b_2$, the interval $[ \log_{a}\left(2 \times 10^P\right),  \log_{a} \left( 3 \times 10^P\right) )$; $\dots$; or, $b_9$, the interval $[ \log_{a} \left(9 \times 10^P\right),  \log_{a}  \left( 10 \times 10^P\right) )$. The main point is this: if $ \log_{a} a^{R_i}$ lies within the bin $b_d$, then the first digit of  $a^{R_i}$ is $d$.

Since $R_i$ was drawn from a uniform distribution, it has the same chance to take any value within the interval $\left[ \left. P\log_a 10, (P +1) \log_a 10 \right) \right.$ and, therefore, the probability of $ \log_a a^{R_i}$ to fall within the bin $b_d$ is the length of the bin $b_d$ divided by the length of the full interval, namely
\begin{eqnarray}
P\left[  \log_a a^{R_i} \in b_d \right] &=& \frac{\log_{a} \left((d + 1) \times 10^P \right) - \log_{a} \left(d \times 10^P \right)}{\log_{a} \left(10^{P+1} \right) - \log_{a} \left(10^P \right)} \nonumber \\
&=& \frac{ \log_{a} \left(\frac{1 + d}{d} \right)}{\log_{a} 10} \nonumber \\
& = & \log_{10} \left( 1 + \frac{1}{d} \right) .\label{BL1}
\end{eqnarray}
This has the form of Benford's Law for the first digit $P_B^{(1)}(d)$, Eq. (\ref{BL}). Thus, as $N \to \infty$, the first digit distribution of the sequence $\{a^{R_1}, a^{R_2}, \dots , a^{R_N}\}$ will approach $P_B^{(1)}(d)$. We shall call this the {\it General Result} (GR). Note that GR is independent of the integer value of $P$  of the interval as long as the sequence is uniformly distributed. Clearly, the result holds if we change the interval to $\left[ \left. P \log_a 10,  (P + M) \log_a 10 \right) \right.$ with $M$ any integer. Note that we never used the fact that neither $a$, nor $R$, nor $d$, are numbers base 10; the graph in Fig. \ref{GR-2} is plotted for numbers base 10 for illustration purposes, but the result would have been the same for any number base. Thus, we conclude that BL is base invariant, i.e. valid for any number base $K$, with $d = 0, 1, 2, \dots, K -1$. 

The second, third, ..., digit distributions follow right away with the same argument. For instance, the probability that the second digit of $a^{R_i}$ is $d$, equals the sum of the lenghts of the ``sub-bins", 
\begin{eqnarray}
b_{1d} & = & \left[\log_{a}\left((1 + \frac{d}{10}) \times 10^P\right), \log_{a} \left((1 + \frac{ d+1}{10}) \times 10^P\right) \right) \nonumber \\
b_{2d} & = &   \left[\log_{a}\left((2 + \frac{d}{10}) \times 10^P\right), \log_{a} \left((2 + \frac{ d+1}{10}) \times 10^P\right) \right)\nonumber \\
& \vdots & \nonumber \\
b_{9d} & = &  \left[\log_{a}\left((9 + \frac{d}{10}) \times 10^P\right), \log_{a} \left((9 + \frac{ d+1}{10}) \times 10^P\right) \right)\nonumber
\end{eqnarray}
where $d$ can now take all values $0, 1, \dots, 9$. Thus, the second digit distribution is,
\begin{eqnarray}
P_B^{(2)}(d) &=& \frac{1}{\log_{a} \left(10^{P+1} \right) - \log_{a} \left(10^P \right)}  \sum_{m=1}^9 \left[  \log_{a} \left((m + \frac{ d+1}{10}) \times 10^P \right) - \log_{a} \left((m + \frac{ d}{10}) \times 10^P \right) \right] \nonumber \\
&=&  \frac{1}{\log_{a} \left(10^{P+1} \right) - \log_{a} \left(10^P \right)} \sum_{m=1}^9   \log_{a}\frac{ 10m + d + 1}{10 m + d }  \nonumber \\
& = & \sum_{m=1}^9  \log_{10} \left( 1 + \frac{ 1}{10 m + d } \right). \label{BL2}
\end{eqnarray}

The argument is easily generalized to the $n$-th digit and the result is,
\begin{equation}
P_B^{(n)}(d) = \sum_{m=10^{n-2}}^{10^{n-1}-1}  \log_{10} \left( 1 + \frac{ 1}{10 m + d } \right). \label{BLn}
\end{equation}

The present derivation is extremely simple. Although Newcomb never wrote the formulas for $P_B^{(n)}(d)$, given his mastery of log tables and numerical analysis\cite{Newcomb-PT}, it is clear that he knew them since he writes the values of $P_B^{(1)}(d)$ and $P_B^{(2)}(d)$ explicitely and mentions how $P_B^{(3)}(d)$ and $P_B^{(4)}(d)$ behave (the latter are almost uniform). Due to Newcomb's most important result, as we discuss in the next section, it seems to the author that he knew a derivation very similar to this one. Benford's Law and its generalizations have been rigorously shown by Hill\cite{Hill} to follow as a consequence of base invariance of the underlying law. We have no pretense of such a mathematical rigour here, but rather to show its simplicity to a wider audience. 

\subsection{Scale invariance.}

A very important property of BL that follows from GR above is the fact that BL is scale invariant\cite{Pietronero}. Add to the values $R_i$ any constant value $c$. This is equivalent to consider a uniform sequence of numbers $R_i$ in the interval $\left[ \left. c +  P \log_a 10,  c + (P +1) \log_a 10 \right) \right.$. Referring to Fig. \ref{GR-2} one can see that in the semi-log graph this also amounts to shift the interval in the ordinate by a constant amount, $a^c$; one also sees, however, that the sizes of the bins $b_n$ remain unchanged. Thus, the sequence $\{a^c a^{R_1}, a^c a^{R_2}, \dots, a^c a^{R_N}\}$ also obeys BL. But this new sequence is the same as the original one $\{a^{R_1}, a^{R_2}, \dots , a^{R_N}\}$ multiplied by a constant, arbitrary, factor $a^c$.

\subsection{The mantissa rule.}

The sequences or sets of numbers that usually follow BL are not of the form $\{ a^R \}$. So, one can enquiry for a rule that tells us if some sequences do follow BL or not. The answer was also given by Newcomb in his two-page paper\cite{Newcomb}. We shall demonstrate now that for  a given sequence of numbers $\{A_1, A_2, \dots , A_N\}$, if the mantissas of the logarithms of $A_i$, namely of $\{ \log_{10} A_1, \log_{10} A_2, \dots , \log_{10} A_N \}$, are uniformly distributed, then the sequence $\{A_1, A_2, \dots , A_N\}$ obeys BL.  Before we give the demonstration, we note that GR can be restated much simpler for the case $a = 10$, namely, if a sequence of numbers $R_i$ are uniformly distributed in the interval $[0, 1)$, the sequence $10^{R_i}$ follows BL. We use this form below.

The demonstration can be done writing the log of $A_i$ as,
\begin{equation}
\log_{10} A_i =  C(A_i) + m(A_i) ,
\end{equation}
where $C(A_i)$ is an integer, the so-called {\it characteristic} of the log, and $m(A_i)$ the fractional part of the logarithm or {\it mantissa}. Note that by definition the mantissas of logarithms base 10 are within the interval $[0, 1)$. It is clear, then, that when taking the ``antilogarithm" $10^{\log_{10} A_i} = 10^{C(A_i)} 10^{m(A_i)}$ the digits of $A_i$ will be determined only by $10^{m(A_i)}$ since the factor $10^{C(A_i)}$ just determines the position of the decimal point. Thus, the distribution of digits is determined by considering the sequence of the mantissas only,  namely of the sequence $\{ m( A_1),m(A_2), \dots , m(A_N )\}$. Hence, if the latter are uniformly distributed, by GR the sequence $\{10^{m(A_1)}, 10^{m(A_2)}, \dots, 10^{m(A_N)} \}$  obeys BL and, therefore, so does the sequence $\{A_1, A_2, \dots , A_N\}$.

This result is very useful since allows us to check if a sequence of numbers obeys BL by looking at the distribution of the mantissas of their logarithms. This is a simple operational rule, instead of a logical one by checking at the digits themselves.

\section{Some sequences that obey BL.}

The next question is which type of sequences or sets of numbers follow BL.  Answering this in an exhaustive fashion appears as a difficult task. Here, we discuss two general type of sequences that can be shown quite clearly that obey BL. With these two, we shall conjecture about the general case. We discuss these cases below.

\subsection{Products of variables with arbitrary distributions.}

Consider the set of numbers $\{Q_1, Q_2, \dots , Q_N\}$, with $Q_i$ given by the product of $M$ numbers,
\begin{equation}
Q_i = R_1^{(i)} R_2^{(i)}  \cdots R_M^{(i)} ,\label{prod}
\end{equation}
where  $R$ are the absolute values of numbers drawn from an {\it arbitrary} distribution (up to a requirement to be given below). We now show that in the limit $N \to \infty$ and $M \to \infty$, the sequence $\{Q_1, Q_2, \dots , Q_N\}$ obeys BL. 

The idea is to use the mantissa rule. For this, we consider the log base 10 of the numbers $Q_i$,
\begin{equation}
\log_{10} Q_i = \log_{10} R_1^{(i)}  + \log_{10} R_2^{(i)}  + \cdots + \log_{10} R_M^{(i)} .
\end{equation} 
We now introduce the requirement that the {\it distribution of the logarithms} of the numbers $R$ have finite first and second moments. Then, in the limit $M \to \infty$, by the Central Limit Theorem (CLT)\cite{Feller}, the distribution of $\log_{10} Q$ is the normal distribution. That is, the values of $\log_{10} Q$ are distributed as,
\begin{equation}
\rho(\log_{10} Q) = \frac{1}{\sqrt{2 \pi} \sigma} \> e^{- (\log_{10} Q - \log_{10} Q_0)^2/2 \sigma^2 }. \label{gauss}
\end{equation}
Note that this is not the log-normal distribution, but simply the normal distribution for the variable $\log_{10} Q$. The centroid $\log_{10} Q_0 \approx M  c_0$ and $\sigma \approx \sqrt{M} \sigma_0$, where $c_0$ and $\sigma_0$ are the first and second moments of the distribution of the {\it logarithms} of the numbers $R$. This point will be further discussed below. 

\begin{figure}
\begin{center}
 \begin{minipage}[h]{.50\textwidth}
  \includegraphics[width=.99\textwidth]{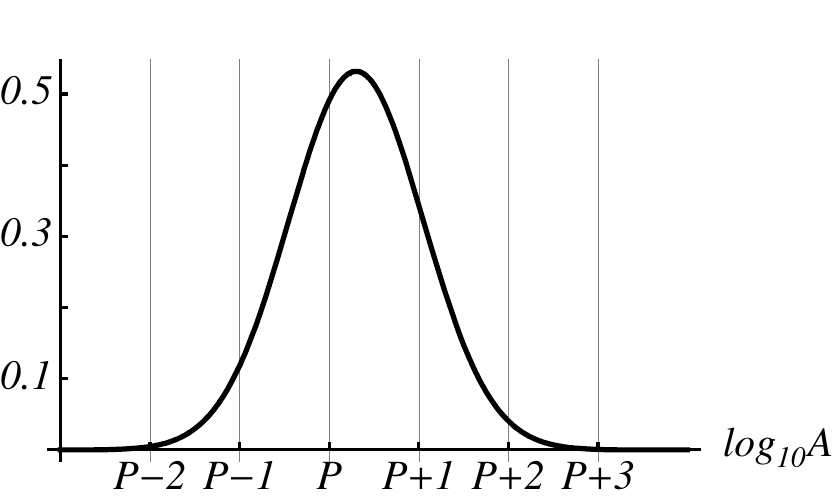}
 \end{minipage}%
 \hfill
 \begin{minipage}[h]{.50\textwidth}
 \centering
  \includegraphics[width=.99\linewidth]{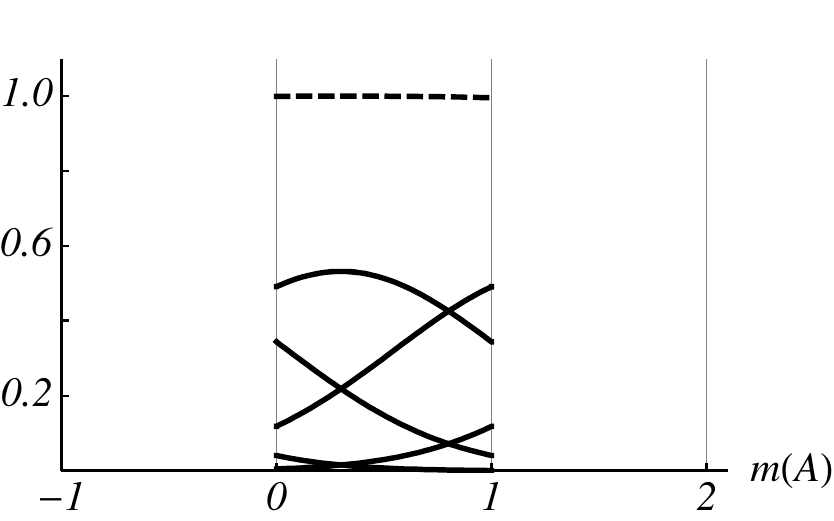}
 \end{minipage}
 \end{center}
\caption{First panel, normal distribution of $\log_{10}(A)$, covering 5 decades approximately. Second panel, the mantissas of the normal distribution within one decade, i.e. in the interval $[0, 1)$; the dotted line is the sum of only 5 decades, adding to 1 within 3 significant figures.}
\label{Gaus}
\end{figure}

We proceed to show that the mantissas of the sequence  $\{\log_{10} Q_1, \log_{10} Q_2, \dots , \log_{10} Q_N\}$ are uniformly distributed in the interval $[0, 1)$, in the limits mentioned. Before we give the general condition, we can see the how this  limit is achieved. Assume that the gaussian function given by Eq.(\ref{gauss}) is already wide enough such that it covers several orders of magnitude, or ``decades'', of the values of $Q$; see Fig. \ref{Gaus} where the decades are denoted by $P - 2$, $P - 1$, $ \dots$, $P + 3$ . The mantissas of the $\log_{10} Q$ are the decimal values within the intervals $P - L$ and $P - (L + 1)$. Thus, we can ``shift'' all intervals within all decades to a single interval, thus placing the mantissas within the same interval. Adding all the values of the mantissas yields, almost, a uniform distribution. This procedure is the same as considering the sum of an infinite number of gaussians each centered at $(\log_{10} Q_0 - P)$ with $P$ taking all the integer values; in the limit $M \to \infty$, equivalent to $\sigma \to \infty$, one gets the exact result,
\begin{equation}
\lim_{\sigma \to \infty} \> \frac{1}{\sqrt{2 \pi} \sigma} \sum_{P= -\infty}^{\infty}  \> e^{- (\log_{10} Q - \log_{10} Q_0 + P)^2/2 \sigma^2 } = 1 .
\end{equation}
This proves that the mantissas are uniformly distributed in the limit, for logarithms normally distributed. Although the previous result is strictly valid only in the limit $\sigma \to \infty$, the convergence is extremely fast. For instance, for $\sigma \approx 1$, the sum differs from 1 in the eighth significant figure. One finds strong deviations from the uniform distribution as $\sigma$ becomes much smaller than 1, that is, when the gaussian covers less than one decade. 

On the other hand, since $\sigma$ depends not only on $M$ but also on the second moment $\sigma_0$ of the distribution of logarithms of $R$, i.e. $\sigma \approx \sqrt{M} \sigma_0$, the convergence might be very slow if the width, or support, of the distribution of $R$ itself is very narrow. As particular examples, considering $R$ taken from a uniform distribution in the interval $[1, 10)$, requires $M$ to be less than 10 (about 4 or 5) to converge to BL. Conversely, for $R$ in the interval $[5, 6)$ takes $M \approx 400$ to yield BL. 

The result of this section, namely of the product of numbers obeying BL, is very robust and general\cite{Pietronero} in the sense that even if the distribution of the numbers $R$ lack second moment, the logarithm of $R$ may not\cite{Lorentz}. This is because the logarithm function is very ``slow'' and tends to smooth the original distribution. Moreover, even if the numbers $R$ are correlated, the action of the logarithm and the limit of very large products (i.e. large values of $M$) may again yield a normal distribution of the logarithms\cite{GUE,Flores}.

\subsection{Generalized geometric sequences of variables uniformly distributed.}

Here we consider a geometric sequence of products of the form, 
\begin{equation}
\{Z^{(1)}, Z^{(2)}, \dots, Z^{(N)}, \dots \} = \{ R_1^{(1)}, R_1^{(2)} R_2^{(2)}, R_1^{(3)}R_2^{(3)}R_3^{(3)}, \dots, 
R_1^{(N)}R_2^{(N)}\cdots R_N^{(N)}, \dots \} ,\label{geo}
\end{equation}
where $R_i^{(J)}$ are numbers uniformly distributed in an arbitrary interval $[a, b]$, with $a$ and $b$ real positive numbers. This sequence obeys BL. Although this result may be generalized to arbitrary distributions, we restrict the results here to uniform distributions. We note that if $a = b$, the above sequence is a true geometric progression with ratio $a$. Thus, geometric progressions also obey BL (except if $a = 10^L$ with $L$ any integer).

\begin{figure}
\begin{center}
 \begin{minipage}[h]{.33\textwidth}
  \includegraphics[width=.99\textwidth]{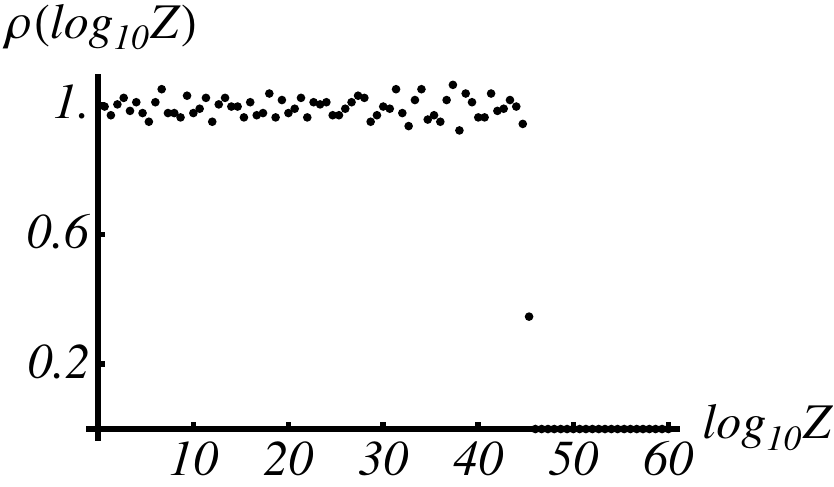}
 \end{minipage}%
 \hfill
 \begin{minipage}[h]{.33\textwidth}
 \centering
  \includegraphics[width=.99\linewidth]{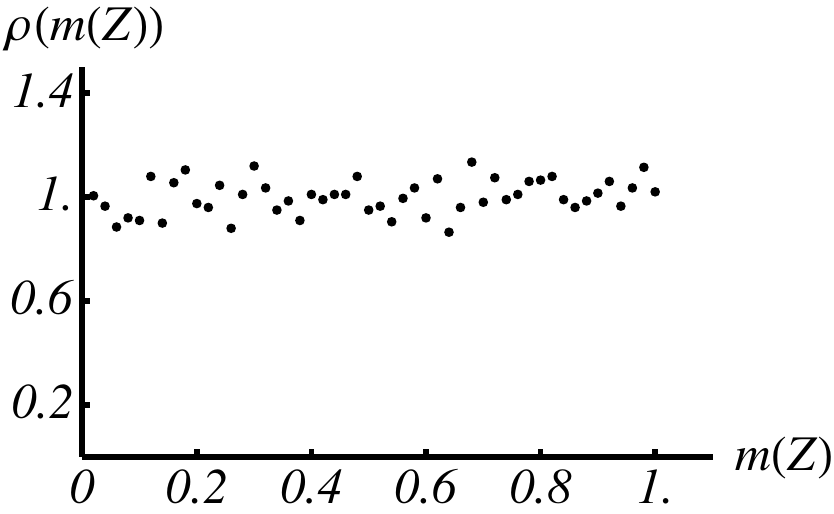}
 \end{minipage}
  \begin{minipage}[h]{.33\textwidth}
  \includegraphics[width=.99\textwidth]{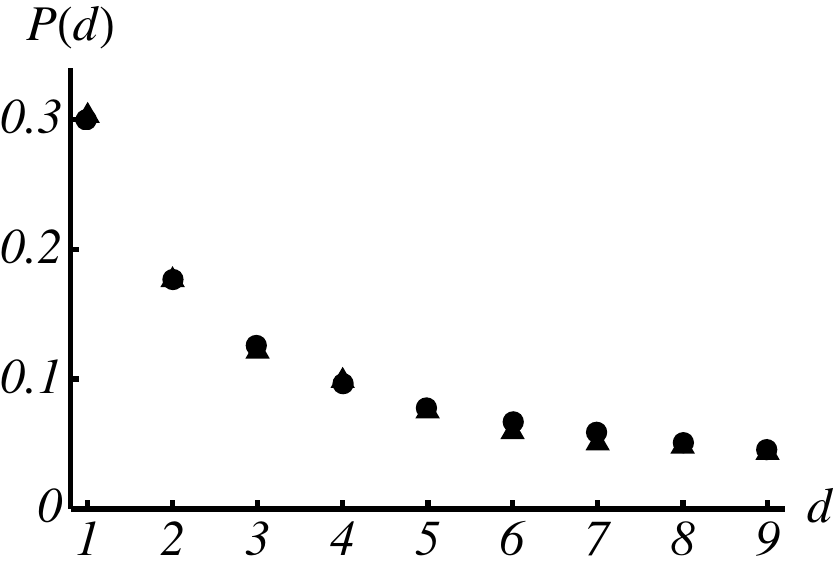}
 \end{minipage}%
 \hfill
 \end{center}
\caption{Numerical analysis of the first 10,000 terms of one realization of a generalized geometric sequence, as given by Eq.(\ref{geo}), for numbers uniformly distributed in the interval $[1.0, 9.9]$. First panel shows the uniform distribution of $\log_{10}(A)$ covering more than 40 decades. Second panel shows that the distribution of mantissas is uniform in $[0, 1)$. The third panel is a comparison of the exact Benford Law (circles) with the distribution of the first digit of the 10,000 terms considered (triangles).}
\label{geofig}
\end{figure}

Again, we first consider the sequence of the logarithms of the products, $\log_{10} Z^{(J)}$. Since we do not have an analytic demonstration, we resort to a numerical one. In Fig. \ref{geofig}  a particular example shows that the distribution of $\log_{10} Z^{(J)}$ becomes uniformly distributed as $J \to \infty$. As this distribution covers many decades of $Z^{(J)}$, obviously the mantissas of $\log_{10} Z^{(J)}$ also become uniform in $[0, 1)$. A numerical comparison with BL is also included. We have extensively verified that these results hold for any sequence of this type, including true geometric progressions\cite{Raimi}.

\subsection{A conjecture on the general case.}

From the above two cases, it appears that a generalization is as follows: As long as the distribution of logarithms is wide enough, namely, covering many decades of the set considered, the mantissa distribution will tend to become uniformly distributed. An analogous argument was recently used by Fewster\cite{Fewster} to illustrate when Benford's law should be obeyed.

\section{Why the pages of log tables wear out following BL?}

Simon Newcomb initiates his article by pointing out that the log tables were worn out more at the beginning than at the end, i.e. following BL.
That is, since the tables are for logarithms of numbers going in order from $1.000\dots$  to $9.999\dots$, he found that the pages for numbers starting with 1 were more used than those for numbers starting with 2, etc. Newcomb gave an explanation of these observation by assuming that ``natural" numbers, i.e. those appearing in Nature, were obtained by ratios of other numbers. Then, he argued that no matter the underlying law of the primitive numbers, their ratios (in the limit of many ratios) had the mantissas of their logarithms uniformly distributed. He then simply stated that this implied Benford Law. As we have seen, the mantissa rule is equivalent to GR. It is fairly evident that Newcomb certainly knew this result, and thus, that he must be credited with the derivation of BL. We mention, once more, that the arguments given in this article are essentially contained in Newcomb's original paper.

An interesting aspect is why Newcomb considered that ``natural" numbers were the result of ratios, or products for that matter, of other numbers. In the light of the previous sections and a bit of second-guessing, we can advance an explanation for this assertion by Newcomb. Moreover, this may also well be the explanation for the agreement of actual real-life data with BL.

To begin, we should recall why log tables were used in the first place. We are well into the era of electronic calculators, be it a pocket-size one or a huge supercomputer: numerical calculations are now their task not ours. But as recently as the early 1970's, not to mention in the XIX century, numerical calculations were done by hand and/or sliding rules. And the log tables were essential to realize those tedious and lengthy tasks. As a matter of fact, logarithms were invented (or discovered?) by John Napier in 1614 to perform lengthy calculations! In the words of Napier himself\cite{e},

{\it Seeing there is nothing that is so troublesome to mathematical practice, nor that doth more molest and hinder calculations, than the multiplications, divisions, square and cubic extractions of great numbers. ... I began therefore to consider in my mind by what certain and ready art I might remove those hindrances.} - John Napier, {\it Mirifici logarithmorum canonis descriptio} (1614).

That is, the trouble appears when one must make calculations by hand, specially multiplications, of numbers with {\it many} digits. It is lengthy, tedious and prone to produce mistakes. Thus, one goes to the tables to find out the logarithms of the numbers involved,  performs sums and subtractions which are much easier, and then taking antilogarithms the result is found. The point is, where did those {\it long} numbers come from? Those were definitely not made up, neither read out from somewhere else, nor measured. The long numbers came themselves from {\it multiplication, divisions or powers} of smaller numbers. The latter may be random, or measured or taken arbitrarily from somewhere else, indeed. But, we insist, the long numbers did arise from operations performed on smaller numbers. As we have seen in the previous section, multiplication of numbers tipically tend to BL, even if only few factors are involved, as long as they arise from wide distributions.  In other words, the numbers that people looked their logarithms for, typically, obeyed already BL. Since in the XIX century numbers were not churned out from a computer but arised from arithmetic operations performed by real people, it seems to the author that for Newcomb these were ``naturally" produced. This may also explain why many sets of real-life data obey BL, that is, unless one asks a computer for a random number,  numbers that quantify a property, be it the area of a lake or the weight of a molecule, usually arise from arithmetic operations performed on measured quantities with arbitrary constants and units. 

{\bf Acknowledgments}. I thank R. Esquivel and A. Robledo for several important references.


\begin{thebibliography}{99}
\bibitem{others} The present article does not pretend to give an exhaustive list of references on BL. The articles by Raimi\cite{Raimi} and Fewster\cite{Fewster} give a good account of several rigorous  and not-so rigorous works. As mentioned in the text, the general result (GR) of Section II appears to be known to many researchers.
\bibitem{Lorentz} One can study the case where $R_i$ is taken from a Lorentzian distribution, that lacks second moment. One finds that the distribution of logarithms of $R_i$ do have a second moment, thus yielding a normal distribution for the logarithms of products of $R_i$.
\end{thebibliography}
\end{document}